\documentstyle[12pt]{article}

\newcommand{\LI}{{\cal{L}}}

\newcommand{\KI}{{\cal{K}}}
\newcommand{\CI}{{\cal{C}}}
\newcommand{\PPI}{{\cal{P}}}
\begin{document}

\title{Hyperbolic cosines and sines theorems for the triangle formed by intersection of three semicircles on Euclidean
plane}
\author{
Robert M. Yamaleev\\
Facultad de Estudios Superiores,\\
Universidad Nacional Autonoma de Mexico\\
Cuautitl\'an Izcalli, Campo 1, C.P.54740, M\'exico.\\
Joint  Institute for Nuclear Research, LIT, Dubna, Russia.\\
 Email:iamaleev@servidor.unam.mx }
 \maketitle
%
%
\begin{abstract}

 The purpose the present paper is to construct the hyperbolic trigonometry on Euclidean plane without refereing to
 hyperbolic plane. In this paper we  show that the concept of hyperbolic angle and its functions forming {\it the hyperbolic
 trigonometry} give arise on Euclidean plane in a natural way. The method is based on a key- formula establishing a relationship
 between exponential function and the ratio of two segments. This formula opens a straightforward pathway to
 hyperbolic trigonometry on the Euclidean plane.  The hyperbolic law of
cosines I and II and the hyperbolic law of sines  are derived by using of the key-formula and the methods of Euclidean
Geometry, only. It is shown that these laws are consequences of the  interrelations between distances and radii of the
 intersecting semi-circles.

\end{abstract}

\section{ Introduction}

Traditionally interrelations between angles and sides of a triangle are described by the trigonometry via periodic
sine-cosine functions. The concept of the angle in Euclidean plane is intimately related with the figure of a circle
and with  motion of a point along the circumference.

The hyperbolic trigonometry also intimately is related with the circle
 reflecting the "hyperbolic" properties of the circles. The properties of a triangle formed by intersection of three
 semi-circles plays a principal role in the {\it upper half-plane} model of the hyperbolic plane ( Lobachevskii plane).
There are several models for the hyperbolic plane \cite{Stahl}, however all these models use the same idea as {\it
upper half-plane model} $H$ \cite{Katok}, \cite{Anderson}.
 The geometry of the hyperbolic plane within the framework of $H$ model is
 studied by considering quantities invariant under an action of the general {\it M\"{o}bius group}, which consists of
 compositions of M\"{o}bius transformations and reflections \cite{Ungar}.
The model of the hyperbolic plane is the {\it upper half-plane} model. The underlining space of this model is the upper
half-plane $H$ in the complex plane $C$ defined to be
$$
H=\{ z\in C:Im(z)>0 \}.
$$
It is used the usual notion of point that $H$ inherits from $C$. It is also used the usual notion of angle that $H$
inherits from $C$, that is, the angle between two curves in $H$ is defined to be the angle between the curves when they
are considered to be curves in $C$, which in turn is defined to be the angle between their tangent lines.
 Usually it is thought that the hyperbolic laws of cosines and sines in the upper half-plane model are
consequences of a special structure of the complex plane.

 The purpose the present paper is to construct the hyperbolic trigonometry on Euclidean plane without refereing to
 hyperbolic plane. The principal  goal of the paper is to  show that the concept of hyperbolic angle and its functions forming {\it the hyperbolic
 trigonometry} give arise on Euclidean plane in a natural way. Our method is based on one simple  but useful formula
 of the  hyperbolic calculus which we denominated as {\it a Key-formula}. The Key-formula establishes relationship
 between exponential function and the ratio of two quantities. This formula opens the straightforward pathway to
 hyperbolic trigonometry on Euclidean plane.  It is well-known, the hyperbolic {\it "law of
sines"} and the hyperbolic {\it "law of cosines II"} are derived from the hyperbolic law of cosines I by algebraic
manipulation \cite{Katok}, \cite{Anderson}. In this paper, we prove all three laws separately on making use of
interrelations between distances between centers and radii of circles on Euclidean plane, only.

In Section 2 the Key-formula of hyperbolic calculus is substantiated. On making use of the Key-formula elements of the
right-angled triangle are expressed via hyperbolic trigonometry.
 In Section 3, two
 main relationships between elements of the triangle formed by intersections of semicircles is established.
 In Section 4,
the theorem of cosines for the triangle formed by three intersecting semicircles, denominated as {\it hyperbolic "law
of cosines I"}, is proved.
 In Section 5,
  the hyperbolic {\it "law of sines"} and the hyperbolic {\it "law of
cosines II"} are derived on making use of the main relationships between thee circles.

\section{ Hyperbolic trigonometry in Euclidean Geometry}

{\bf 2.1 Key-formula of hyperbolic calculus}.

The Key-formula of hyperbolic calculus has been established in Ref.\cite{Yamaleev}.
 Let $a=a(s),~b=b(s)$ are real functions of the parameter $s$. If the difference $(a-b)$ does not depend of the
 parameter $s$ then the following formula holds true
$$
\exp((a-b)s)=\frac{a}{b}. \eqno(2.1)
$$
The formula (2.2) constitutes a substance what we shall refer as {\it Key-formula}.
 The main advantage of the Key-formula is the following: {\it the argument of the exponential function is proportional
to the difference between the nominator and the denominator, if this difference is a constant}.

If $(a-b)$ does not depend of s then functions $a$ and $b$ can be presented as follows
$$
a=a_0+Y(s),~~b=b_0+Y(s),
$$
where $a_0,b_0$ are constants. For the ratio of these constants find $\phi_0$ satisfying the equation
$$
\exp((a_0-b_0)\phi_0)=\frac{a_0}{b_0}.
$$
Since a shift of the pair $\{ a_0,b_0 \}$ by $Y$ does not change $(a_0-b_0)$, one may write,
$$
\frac{a_0+Y}{b_0+Y}=\frac{a}{b}=\exp((a_0-b_0)(\phi_0+\delta))=\exp((a-b)\phi). \eqno(2.2)
$$

{\bf 2.2 Elements of right-angled triangle as functions of a hyperbolic argument }.

Let $\triangle ABC$ be a right-angled triangle with right angle at $C$. Denote the sides by $a,b$, the hypotenuse
 by $c$, the angles opposite by $A$, $B$, $C$, correspondingly. Circular functions of the angles are defined
 in the usual way. For instance,
$$
\sin B=\frac{b}{c},~~\tan B=\frac{b}{a}. \eqno(2.3)
$$
These two ratios are functions of the same angle $B$. Our purpose is to introduce an hyperbolic angle in a such way
that these ratios will be expressed as functions of unique hyperbolic angle. This achieved as follows.  Define the
hyperbolic angle $\xi$ by the following relation
$$
\frac{c+a}{c-a}=\exp(2\xi).    \eqno(2.4)
$$
From this equation it follows
$$
\frac{a}{c}=\tanh(\xi),~~\frac{a}{b}={\sinh(\xi)}.   \eqno(2.5)
$$
This leads to the following interrelations between circular and hyperbolic trigonometry
$$
\cos B=\frac{a}{c}=\tanh(\xi),~~\tan A=\frac{b}{a}=\frac{1}{\sinh(\xi)}.   \eqno(2.6)
$$

Now let us explore the same problem by using a geometrical motion. We shall change sides $b$ and $c$  remaining
unchanged the side $a$ and the right angle $C$. Since the length $a$ is a constant of this motion, in agreement with
Key-formula (2.2), we write
$$
\frac{c+a}{c-a}=\exp(2a\phi).   \eqno(2.7)
$$
Notice, now the argument of exponential function in (2.6) is proportional to $a$: $\xi=a\phi$, where $\phi$ is a
parameter of the evolution. Formulae (2.6) are re-written as
$$
\cos B=\frac{a}{c}=\tanh(a\phi),~~\tan B=\frac{b}{a}=\frac{1}{\sinh(a\phi)}.   \eqno(2.8)
$$

It is interesting to observe that within these representation we are able to find a limit for $a=0$. In fact, at this
limit we get
$$
c(a=0)=b(a=0)=\frac{1}{\phi}.
$$
The relations between circular and hyperbolic trigonometry (2.8) also may be put into other form, for instance,
$$
~\cot B=\sinh(a\phi),  \eqno(2.9)
$$
or,
$$
\tan \frac{B}{2}=\exp(-a\phi). \eqno(2.10)
$$

Install the triangle $\triangle ABC$ in such a way that the side $b$ lies on horizontal axis $X$, the side $a$ is
perpendicular to this line at the point $C$ (Fig.1). The line $c$ cuts $X$ in $A$. The point of intersection $A$ may
continue moving, and the distance $AC$ tends to infinity.

Let us recall the problem of parallel lines in Geometry.

 The ray $AB$ then tends to a definite limiting
position $BL$, and $BL$ is said to be parallel to $X$. As the point $A$ moves along $X$ away from $C$ there are two
possibilities to consider:

(1) In Euclidean Geometry, the angle between two lines $BL$ and $BC$ is equal to right angle.

(2) The hypothesis of Hyperbolic Geometry is that this angle less than the right angle.

The most fundamental formula of the Hyperbolic Geometry is the formula connecting the parallel angle $\Pi(d)$ and the
length $d$. In order to establish those relationship the concept of horocycles, some circles with center and axis at
infinity, were introduced \cite{Lobatschewsky}. The great theorem which enables one to introduce the circular
functions, sines and cosines of an angle is that the geometry of shortest lines (horocycles) traced on horosphere is
the same as plane Euclidean geometry. The function connecting the parallel angle with the distance $d$ is given by
$$
\exp(-\frac{d}{\kappa})=\tan\frac{\Pi(p)}{2}. \eqno(2.11)
$$
Now, introduce in (2.10) the value inverse to $\phi$, $ \KI=\frac{1}{\phi},$ and write (2.10) in the form
$$
\exp(-\frac{a}{\KI})=\tan\frac{B}{2}. \eqno(2.12)
$$

 Let the point to tend to infinity. The following  two cases can be considered.

(1) $\phi$ will tend to zero, $\KI\rightarrow \infty$, and $\angle B \rightarrow  \frac{\pi}{2}$.

This is true in Euclidean Geometry.

 (2) Suppose that $\phi$, $\KI$ and the angle $B$ go to some limited  values,
$$
 \lim_{AC\rightarrow \infty}\KI=\kappa,~~\lim B=\Pi(a). \eqno(2.13)
$$
In this way, from (2.10) we come to main formula of hyperbolic geometry.

{\bf 3.2 Rotational motion of a line tangent to the semicircle.}

Let us explore a motion of the line $\LI$ around the semicircle ${\CI}$  remaining tangent to the semicircle
\cite{Yamaleev}.

Draw semicircle  $\CI$  (Fig.2) end-points and the center of which lie on horizontal axis, $X$-axis. Denote by $B$
center of the semicircle and by $K_1$, $K_2$ end-points of the semicircle on $X$-axis. Draw the line $\LI$ tangent to
the semicircle at the point $C$. This line intersects with X-axis at the point $A$. Through end-points of the
semicircle, $K_1$, $K_2$, erect the lines parallel to vertical axis, $Y$-axis. The intersections of these lines with
line $\LI$ denote by $P_1$ and $P_2$, correspondingly. Draw line parallel to Y-axis from the center $B$ which acroses
$\CI$ at the point $N$, at the top of the semicircle $\CI$. This line intersects with the line $\LI$ at the point
$M_1$. From point $M_1$ draw horizontal line which acroses the vertical line $K_2P_2$ at the point $M_2$.

Denote by $r$ radius of the circle, so that $r=BN=BC$. Denote by $\angle B$ the angle between $AB$ and $BC$. The angle
$\angle C$ is rectangle, so that
$$
(AB)^2-(AC)^2=r^2. \eqno(2.14)
$$
Consider rotational motion of the line $\LI$ around semicircle $\CI$ remaining tangent to $\CI$ where the point $C$
runs between points $K_1$ and $N$.
 During of the motion the length of segments $AC$ and $AB$ will change, but the triangle remains to be right-angled.
This is exactly the case considered above, difference is that now $\triangle ABC$ installed in another position with
respect to horizontal $X$-axis. We have seen, this evolution process is described by equations
$$
AB=r\coth(\xi),~~AC=\frac{r}{\sinh(\xi)}, \eqno(2.15)
$$
where $\xi=r\phi$.

From similarity of  triangles $\triangle ABC$ and $\triangle AK_1P_1$ we find
$$
\frac{P_1K_1}{AP_1}=\frac{r}{AC}=\sinh(\xi).\eqno(2.16)
$$
Consequently,
$$
P_1K_1=AK_1\sinh(\xi)=(AB-r)\sinh(\xi)=r\exp(-\xi).\eqno(2.17)
$$
Then, obviously,
$$
P_2K_2=r\exp(\xi).\eqno(2.18)
$$
Furthermore,
$$
BM_1=\frac{P_1K_1+P_2K_2}{2}=r\cosh(\xi),~~P_2M_2=r\sinh(\xi), \eqno(2.19)
$$
Thus, all segments of the lines in Fig.2 can be expressed via hyperbolic angle $\xi$ and the radius of the semicircle
$r$. On the basis of obtained formulae the following relationships between circular and hyperbolic trigonometric
functions are established
$$
\sin(B)=\frac{1}{\cosh(\xi)},~~~\cot(B)={\sinh(\xi)}~~~\cos(B)={\tanh(\xi)}.\eqno(2.20)
$$

\section{ Relationships between elements of three intersecting semicircles}

{\bf 3.1 Hyperbolic cosine- sine functions of arcs of the semicircle.}

Up till now we were able to define hyperbolic trigonometric functions of the arcs originated from the top $N$ ( for
variable $\xi$).
 Now let us calculate trigonometric functions of the arcs with arbitrary end-points on the semicircle.
Consider arc $\widehat{A_1A_2}$ defined in the first quadrant of the semicircle with end-points at the points $A_1$ and
$A_2$ where $\widehat{NA_2}< \widehat{NA_1}$.  Since we are able to calculate arcs with the origin installed on the top
of the semicircle ( for variable $\xi$), we shall present this segment on the circle as difference of two segments both
originated from the top of the semicircle. For example,
$$
\widehat{A_1A_2}=\widehat{NA_1}-\widehat{NA_2}.
$$
Then,
$$
\cosh\widehat{A_1A_2}= \cosh \widehat{NA_1}\cosh \widehat{NA_2} -\sinh\widehat{NA_1} \sinh\widehat{NA_2},
$$
$$
\sinh\widehat{A_1A_2}=\sinh\widehat{NA_1}\cosh \widehat{NA_2} -\cosh\widehat{NA_1}\sinh\widehat{NA_2}.\eqno(3.1)
$$
Denote by $a_1,a_2$ the angles formed by radiuses $O_aA_1$ and $O_aA_2$ with $X$-axis,correspondingly. Then the
functions
$$
\cosh \widehat{NA_1},~\cosh \widehat{NA_2},~\sinh\widehat{NA_1},~ \sinh\widehat{NA_2},
$$
are expressed via circular trigonometric functions as follows
$$
\cosh \widehat{NA_1}=\frac{1}{\sin a_1},~~\cosh\widehat{NA_2}=\frac{1}{\sin a_2},~~\sinh\widehat{NA_1}=\cot a_1,~~
\sinh \widehat{NA_2}=\cot a_2.\eqno(3.2)
$$
On making use of equations (3.2) in (3.1), we get
$$
\cosh \widehat{A_1A_2}= \frac{1-\cos a_1\cos a_2}{\sin a_1\sin a_2},~~~ \sinh \widehat{A_1A_2}= \frac{\cos a_1 -\cos
a_2}{\sin _1\sin a_2}.\eqno(3.3)
$$

{\bf 3.2 The main relationships between elements of the triangle formed by intersections of semicircles.}

In Fig.3 three intersecting semicircles with centers installed on horizontal axis at the points $O_a,O_b,O_c$ are
presented. Intersections of the circumferences form triangle $\widetilde{\triangle} ABC$ bounded by segments of the
circumferences $a=\widehat{BC},~c=\widehat{AB},b=\widehat{AC}$. Connect vertex of the triangle with centers of the
circle by corresponding radiuses. Denote by $a_k,b_k,c_k,k=1,2$ the angles bounded by the radiuses  and $X$-axis, where
$a_1>a_2>0,b_1>b_2>0,c_1>c_2>0$. By making use of (3.3) define hyperbolic cosine-sine functions corresponding to
bounding segments
$$
\cosh a=\frac{1-\cos a_1\cos a_2}{\sin a_1\sin a_2},~~\sinh a=\frac{\cos a_1-\cos a_2}{\sin a_1\sin a_2},\eqno(3.4a)
$$
$$
\cosh b=\frac{1-\cos b_1\cos b_2}{\sin b_1\sin b_2},~~\sinh b=\frac{\cos b_1-\cos b_2}{\sin b_1\sin b_2},\eqno(3.4b)
$$
$$
\cosh c=\frac{1-\cos c_1\cos c_2}{\sin c_1\sin c_2},~~\sinh c=\frac{\cos c_1-\cos c_2}{\sin c_1\sin c_2}.  \eqno(3.4c)
$$
It is used the usual notion of the angle, that is, the angle between two curves  is defined as an angle between  their
tangent lines. Let the angles $\alpha,\beta,\gamma$ be angles at the vertex $A,B,C$, correspondingly. For these angles
we can define its proper cosine and sine functions. The angles of the triangle $\widetilde{\triangle} ABC$
$\alpha,\beta,\gamma$, are closely related with angles $a_1,a_2,b_1,b_2,c_1,c_2$. From draught in Fig.3 we find the
following relationships between them
$$
\beta=a_2+c_2,~\delta=\pi-a_1-b_2,~~\alpha=b_1-c_1.~\eqno(3.5)
$$
Then,
$$
\cos \alpha=\cos b_1\cos c_1+\sin b_1\sin c_1,\eqno(3.6a)
$$
$$
\cos \beta=\cos a_2\cos c_2-\sin a_2\sin c_2,\eqno(3.6b)
$$
$$
\cos \delta=-\cos b_2\cos a_1+\sin b_2\sin a_1,\eqno(3.6c)
$$
$$
\sin \alpha=\sin b_1\cos c_1-\cos b_1\sin c_1,\eqno(3.7a)
$$
$$
\sin \beta=\sin a_2\cos c_2+\cos a_2\sin c_2,\eqno(3.7b)
$$
$$
\sin \delta=\sin b_2\cos a_1+\cos b_2\sin a_1.\eqno(3.7c)
$$

Denote distances between centers by
$$
O_{cb}=O_cO_b,~O_{ba}=O_bO_a,~O_{ac}=O_aO_c.
$$
The theorem of sines employed for triangles $O_cAO_b$, $O_bCO_a$, $O_cBO_a$ gives six relations of type
$$
\frac{\sin\alpha}{O_{cb}}=\frac{\sin c_1}{r_b}=\frac{\sin b_1}{r_c},\eqno(3.8a)
$$
$$
\frac{\sin\gamma}{O_{ba}}=\frac{\sin a_1}{r_b}=\frac{\sin b_2}{r_a},\eqno(3.8b)
$$
$$
\frac{\sin\beta}{O_{ac}}=\frac{\sin c_2}{r_a}=\frac{\sin a_2}{r_c}.\eqno(3.8c)
$$
From these relations  it follows the first set of main relationships: \\
 {\bf  Relation I}.
$$
r_a\sin a_1=r_b\sin b_2,~r_c\sin c_2=r_a\sin a_2,~r_c\sin c_1=r_b\sin b_1,~~\eqno(3.9)
$$
From the draught in Fig.3 it is seen that
$$
O_{ac}=O_{ba}+O_{cb},\eqno(3.10)
$$
where
$$
O_{ac}=r_c\cos c_2+r_a\cos a_2,~~~O_{ba}=r_a\cos a_1+r_b\cos b_2.\eqno(3.11)
$$
Hence,
$$
O_{cb}=r_c\cos c_2+r_a\cos a_2-r_a\cos a_1-r_b\cos b_2.\eqno(3.12)
$$
From vertices of $\widetilde{\triangle} ABC$  erect lines perpendicular to horizontal line, which intersect $X$-axis at
points $ h_A,h_B,h_C $, correspondingly.
 From the draught in Fig.3 we find that
$$
O_{cb}=O_ch_A-h_AO_b=r_c\cos c_1-r_b\cos b_1.\eqno(3.13)
$$
By equating (3.12) with (3.13) we arrive to another main  relationship between radii and angles:\\

{\bf Relation II}.
$$
r_c\cos c_1-r_b\cos b_1= r_c\cos c_2+r_a\cos a_2-r_a\cos a_1-r_b\cos b_2.\eqno(3.14)
$$

We shall effect a simplification by using the following designations.
$$
r_a\cos a_1=w_{01},~r_a\sin a_1=w_1,~~r_a\cos a_2=v_{01},~r_a\sin a_2=v_1,
$$
$$
r_b\cos b_1=w_{02},~r_b\sin b_1=w_2,~~r_b\cos b_2=v_{02},~r_b\sin b_2=v_2,
$$
$$
r_c\cos c_1=w_{03},~r_c\sin c_1=w_3,~~r_c\cos c_2=v_{03},~r_c\sin c_2=v_3.
$$
In these designations  formulae (3.6a,b,c) and (3.7a,b,c) are written as follows
$$
r_br_c\sin \alpha=w_2w_{03}-w_{02}w_3,~~~ r_br_c\cos \alpha=w_{02}w_{03}+w_{2}w_3, \eqno(3.15a)
$$
$$
r_ar_c\sin \beta=v_1v_{03}+v_{01}v_3,~~~ r_ar_c\cos \beta=v_{01}v_{03}-v_{1}v_3,  \eqno(3.15b)
$$
$$
r_ar_b\sin \delta=v_2w_{01}+v_{02}w_1,~~r_ar_b\cos \delta=v_{2}w_{1}-v_{02}w_{01},\eqno(3.15c)
$$
Formulae (3.4a,b,c) for hyperbolic sines and cosines are re-written as follows
$$
\cosh a=\frac{r_a^2-w_{01}v_{01}}{w_1v_1},~~~ \sinh a=r_a\frac{w_{01}-v_{01}}{w_1v_1},
 $$
 $$
\cosh b=\frac{r_b^2-w_{02}v_{02}}{w_2v_2},~~~ \sinh b=r_b\frac{w_{02}-v_{02}}{w_2v_2},
 $$
 $$
\cosh c=\frac{r_c^2-w_{03}v_{03}}{w_3v_3},~~~ \sinh c=r_c\frac{w_{03}-v_{03}}{w_3v_3}.\eqno(3.16)
$$

The equations of main {\bf Relation I} now take the form
$$
x:=w_2=w_3,~y:=v_1=v_3,~z:=w_1=v_2. \eqno(3.17)
$$
Equations (3.8)-(3.11) are re-written as follows
$$
O_{cb}=w_{03}-w_{02},~ O_{ac}=v_{03}+v_{01},~ O_{ba}=w_{01}+v_{02}.
$$
Correspondingly, the main {\bf Relation II}  takes the form
$$
w_{03}-w_{02}=v_{03}+v_{01}-w_{01}-v_{02}.\eqno(3.18)
$$
This expresses the fact that $O_{ca}$ is a sum of $O_{cb}$ and $O_{ba}$. Notice, equation (3.18)  can be re-written
also in another equivalent form, namely,
$$
w_{03}-v_{03}=w_{02}-v_{02}-(w_{01}-v_{01}).\eqno(3.19)
$$
Denote the segments -projections of sides of $\widetilde{\triangle} ABC$  on $X$-axis by $\PPI(AC)= h_Ah_C$,
$\PPI(AB)=h_Ah_B$, $\PPI(BC)=h_Bh_C$. From the draught in Fig.3 it is seen that
$$
\PPI(AC)=\PPI(AB)+\PPI(BC), \eqno(3.20)
$$
where
$$
\PPI(BC)=w_{01}-v_{01},~\PPI(AC)=w_{02}-v_{02},~\PPI(AB)=w_{03}-v_{03}.
$$

\section{Hyperbolic law of cosines I for the triangle formed by intersection of three semicircles }

 The main aim of this section is to prove the hyperbolic law - theorem of cosines I for  triangle
 $\widetilde{\triangle} ABC$ formed by intersection of semicircles with centers installed on $X$-axis (Fig.3) which is
 given by the set of three equations
$$
\cosh c=\cosh a\cosh b-\sinh a \sinh b \cos\delta,
$$
$$
\cosh b=\cosh a\cosh c-\sinh c \sinh a \cos\beta,
$$
$$
\cosh a=\cosh c\cosh b-\sinh c \sinh b \cos\alpha.
$$

{\bf Theorem of cosines I.}

{\it The following equation for elements of the triangle $\widetilde{\triangle} ABC$ formed by intersection of three
circles holds true}
$$
\cosh c=\cosh a\cosh b-\sinh a \sinh b \cos\delta.\eqno(4.1)
$$

{\bf Proof}

 Square both sides of the main  {\bf Relation II} to obtain
$$
(w_{03}-v_{03})^2=(w_{02}-v_{02})^2+(w_{01}-v_{01})^2- 2(w_{02}-v_{02})(w_{01}-v_{01}),  \eqno(4.2)
$$
and evaluate this equality by taking into account formulae (3.15)-(3.16). First of all, evaluate the left-hand side of
this as follows
$$
v^2_{03}+w^2_{03}-2v_{03}w_{03}=~ v^2_{03}+w^2_{03}-2r_c^2+2(r_c^2 -v_{03}w_{03}), \eqno(4.3)
$$
and, notice that
$$
2r_c^2= w_{03}^2+w^2_{3}+v^2_{03}+v^2_{3}.\eqno(4.4)
$$
 Transform (4.3) into the following form
$$
v^2_{03}+w^2_{03}-2r_c^2=v^2_{03}+w^2_{03}-(w_{03}^2+w^2_{3}+v^2_{03}+v^2_{3})=-(w^2_{3}+v^2_{3}).   \eqno(4.5)
$$
 Equation (4.2) is written as
$$
-\underbrace{(w^2_{3}+v^2_{3})} +2(r_c^2 -v_{03}w_{03})~=~(w_{02}-v_{02})^2+(w_{01}-v_{01})^2-
2(w_{02}-v_{02})(w_{01}-v_{01}).          \eqno(4.6)
$$
The underlined term pass from the left-hand side  to the right-hand side of the equation. Then in the left-hand side we
remain with the expression
$$
2(r_c^2 -v_{03}w_{03})=2~v_3w_3~\cosh c.        \eqno(4.7)
$$
We arrive to the following equation
$$
2~v_3w_3\cosh c= \frac{1}{v_3w_3}(~ w^2_{3}+v^2_{3}+(v_{02}-w_{02})^2+(v_{01}-w_{01})^2-
2(v_{02}-w_{02})(v_{01}-w_{01})~) .         \eqno(4.8)
$$
 From the second main {\bf Relation I} we use $w_3=w_2,~~v_3=v_1$. This makes true the following equation
$$
\frac{1}{w_3v_3}=\frac{1}{w_2v_2}\frac{1}{w_1v_1}v_2w_1. \eqno(4.9)
$$
The factor of  the right-hand side of (4.8) replace by the right-hand side of (4.9). In this way we come to the
following equation
$$
2\cosh c=
$$
$$\frac{1}{w_2v_2}\frac{1}{w_1v_1}\{~  v_2w_1(~w^2_{3}+v^2_{3}+(v_{02}-w_{02})^2+(v_{01}-w_{01})^2~) \}-
\underbrace{2\frac{1}{w_2v_2}\frac{1}{w_1v_1}~v_2w_1~(v_{02}-w_{02})(v_{01}-w_{01})}.   \eqno(4.10)
$$
Evaluate now underlined term in the right-hand side of equation (4.10), which we firstly transform as follows
$$
2(v_{02}-w_{02})(v_{01}-w_{01})\frac{1}{w_2v_2}\frac{1}{w_1v_1}v_2w_1= 2 \frac{v_2w_1}{r_ar_b}~
\frac{r_b(v_{02}-w_{02})}{w_2v_2}\frac{r_a(v_{01}-w_{01})}{w_1v_1}. \eqno(4.11)
$$
Then, use the second of equations of (3.15c) written as
$$
\frac{v_2w_1}{r_ar_b}=\cos \delta+\frac{v_{02}w_{01}}{r_ar_b}.         \eqno(4.12)
$$
 By making use of equation (4.12) we are able to evaluate equation (4.11) as follows
$$
2 \frac{v_2w_1}{r_ar_b}~ \frac{r_b(v_{02}-w_{02})}{w_2v_2}\frac{r_a(v_{01}-w_{01})}{w_1v_1}=
$$
$$
2(~\cos \delta+\frac{v_{02}w_{01}}{r_ar_b})~ \frac{r_b(v_{02}-w_{02})}{w_2v_2}\frac{r_a(v_{01}-w_{01})}{w_1v_1}=
$$
$$
2\cos \delta~ \frac{r_b(v_{02}-w_{02})}{w_2v_2}\frac{r_a(v_{01}-w_{01})}{w_1v_1}-
2v_{02}w_{01}\frac{(v_{02}-w_{02})}{w_2v_2}\frac{(v_{01}-w_{01})}{w_1v_1}=
$$
$$
=\underbrace{2\cos \delta~ \sinh a \sinh b}~-
2v_{02}w_{01}\frac{(v_{02}-w_{02})}{w_2v_2}\frac{(v_{01}-w_{01})}{w_1v_1}.\eqno(4.13)
$$
Replace underlined term of (4.10) by (4.13), and  pass the underlined expression of (4.13) to the left-hand side of
obtained equation. As a result, we come to the following equation
$$
2\cosh c-2\cos \delta~ \sinh b\sinh a=
$$
$$
\frac{1}{w_2v_2w_1v_1}~\{~v_2w_1~(~w^2_{2}+v^2_{1}+(w_{02}-v_{02})^2+(w_{01}-v_{01})^2)-
2v_{02}w_{01}(w_{02}-v_{02})(w_{01}-v_{01})~\}.   \eqno(5.14)
$$
On making use of elementary algebra one may show that (see, the section Appendix),
$$
{v_2w_1}~(~w^2_{2}+v^2_{1}+(w_{02}-v_{02})^2+(w_{01}-v_{01})^2~)-2v_{02}w_{01}(w_{02}-v_{02})(w_{01}-v_{01})=
$$
$$
2(r_a^2-v_{01}w_{01})(r_b^2-v_{02}w_{02}).    \eqno(4.15)
$$
 Now,  substitute (4.15) into (4.14) and take into account (3.16). This gives
$$
2\frac{1}{w_1v_1v_2w_2}~(r_a^2-v_{01}w_{01})(r_b^2-v_{02}w_{02})= 2\cosh a~\cosh b,        \eqno(4.16)
$$
by using of which we arrive from  (4.14) to the following equation
$$
2\cosh c-2\cosh a\cosh b=2\cos\delta\sinh a\sinh b. \eqno(4.17)
$$

{\bf End of Proof}.

The other two equations of the law, obviously, are proved analogously.

\section{Hyperbolic laws of sines and cosines II}

The main task of this section is to prove {\it Hyperbolic law (theorem) of sines }, which is given by the formulae
$$
\frac{\sinh a}{\sin\alpha}= \frac{\sinh b}{\sin\beta}=\frac{\sinh c}{\sin\delta},\eqno(5.1)
$$
and the {\it Hyperbolic law (theorem) of cosines II} given by the formulae
$$
\cos\delta=-\cos\alpha\cos\beta-\sin\alpha \sin\beta \cosh c,\eqno(5.2)
$$
$$
\cos\beta=-\cos\alpha\cos\delta-\cosh b \sin\alpha \sin\delta, \eqno(5.3)
$$
$$
\cos\alpha=-\cos\beta\cos\delta-\cosh a \sin\delta\sin\beta.\eqno(5.4)
$$

{\bf 5.1 Hyperbolic theorem of sines and its geometrical interpretation on Euclidean plane}.

{\bf Lemma 5.1}

{\it The ratios of projections of the sides of triangle $\widetilde{\triangle} ABC$ on $X$-axis to corresponding
distances between centers of the semicircles are equal to each other.}

{\bf Proof}

Projections of the sides of $\widetilde{\triangle} ABC$ are given by formulae
$$
\PPI(BC)=r_a\cos a_1-r_a\cos a_2,~\PPI(AC)=r_b\cos b_1-r_b\cos b_2,~\PPI(AB)=r_c\cos c_1-r_c\cos c_2.\eqno(5.5)
$$
Distances between centers of the circles have been defined as (see, (3.11), (3.12)),
$$
O_{ca}=r_c\cos c_2+r_a\cos a_2,~~~O_{ba}=r_a\cos a_1+r_b\cos b_2,~~~ O_{cb}=r_c\cos c_1-r_b\cos b_1,\eqno(5.6)
$$
and,
$$
{\PPI(AC) }=\PPI(BC)+\PPI(AB),~~   O_{ca}= {O_{bc}+O_{ab}}.\eqno(5.7)
$$
Write the first main  Relation I given by the set of equations
$$
r_a\sin a_1=r_b\sin b_2,~r_c\sin c_2=r_a\sin a_2,~r_c\sin c_1=r_b\sin b_1,\eqno(5.8)
$$
 in a squared form, namely,
$$
r^2_a-r^2_a\cos^2 a_1=r^2_b-r^2_b\cos^2 b_2,~r^2_c-r^2_c\cos^2 c_2=r^2_a-r^2_a\cos^2 a_2,~r^2_c-r^2_c\cos^2
c_1=r^2_b-r^2_b\cos^2 b_1. \eqno(5.9)
$$
Then, for the squared distances $O^2_{ik},i,k=a,b,c$ we have,
$$~~~O^2_{ca}=r^2_c\cos^2 c_2+r^2_a\cos^2 a_2+2r_c\cos c_2r_a\cos a_2,
$$
 $$~~~O^2_{ba}=r^2_a\cos^2 a_1+r^2_b\cos^2 b_2+2r_a\cos a_1r_b\cos b_2,
 $$
 $$~~~O^2_{cb}=r^2_c\cos^2 c_1+r^2_b\cos^2 b_1-2r_c\cos c_1r_b\cos b_1. \eqno(5.10)
 $$
 Combine  equations (5.9) with (5.10), this leads to the following system of equations
$$
(a)~~~O^2_{ca}=r^2_c-r^2_a+2r_a\cos a_2~O_{ac},~~~(b)~~O^2_{ca}=r^2_a-r^2_c+2r_c\cos c_2~O_{ac},\eqno(5.11a)
$$
 $$(a)~~~O^2_{ba}=r^2_a-r^2_b+2r_b\cos b_2~O_{ba},~~~(b)~~O^2_{ba}=r^2_b-r^2_a+2r_a\cos a_1~O_{ba},\eqno(5.11b)
 $$
 $$(a)~~~O^2_{cb}=r^2_b-r^2_c+2r_c\cos c_1O_{cb},~~~(b)~~O^2_{cb}=r^2_c-r^2_b-2r_b\cos b_1~O_{cb}.\eqno(5.11c)
 $$
From these equations the cosines of the angles $a_1,a_2,b_1,b_2,c_1,c_2$ are expressed:
$$
\frac{O^2_{ca}-r^2_c+r^2_a}{2r_aO_{ac}}=\cos a_2~,~~~\frac{O^2_{ca}-r^2_a+r^2_c}{2r_c~O_{ac}} =\cos c_2,
$$
 $$
 \frac{O^2_{ba}-r^2_a+r^2_b}{2r_b~O_{ba}} =\cos b_2,~~~\frac{O^2_{ba}-r^2_b+r^2_a}{2r_a~O_{ba}}=\cos a_1,
 $$
 $$
 \frac{O^2_{cb}-r^2_b+r^2_c}{2r_cO_{cb}}= \cos c_1,~~~\frac{-O^2_{cb}+r^2_c-r^2_b}{2r_b~O_{cb}} =\cos  b_1.    \eqno(5.12)
 $$
 Having these formulae  we may present the projection $ {\PPI_{ca}}$ as follows
$$
~\PPI_{ac}=r_b\cos b_1-r_b\cos b_2=
$$
$$
\frac{-O^2_{cb}+r^2_c-r^2_b}{2~O_{cb}}-\frac{O^2_{ba}-r^2_a+r^2_b}{2~O_{ba}}=
$$
$$
\frac{(~-O^2_{cb}+r^2_c-r^2_b~)O_{ba}-(O^2_{ba}-r^2_a+r^2_b~)O_{cb}}{2~O_{cb}~O_{ba}}=
$$
$$
\frac{~-O^2_{cb}O_{ba}+r^2_cO_{ba}-r^2_b~O_{ba}-~~~(~O^2_{ba}-r^2_a+r^2_b~)~O_{cb}}{2~O_{cb}~O_{ba}}=
$$
$$
\frac{~-O^2_{cb}O_{ba}+r^2_cO_{ba}-r^2_b~O_{ba}-~O^2_{ba}O_{cb}+r^2_aO_{cb}-r^2_b~O_{cb}}{2~O_{cb}~O_{ba}}
$$
The first ratio is presented as follows
$$
\frac{\Pi_{ca}}{O_{ca}}= ~ \frac{~-O_{cb}O_{ba}(~O_{cb}+~O_{ba}) +r^2_cO_{ba}-r^2_b~(O_{ba}-~O_{cb})
+r^2_aO_{cb}}{2~O_{cb}~O_{ba}O_{ca}}=
$$
$$
\frac{~-O_{cb}O_{ba}(~O_{ca}) +r^2_cO_{ba}-r^2_b~(O_{ca}) +r^2_aO_{cb}}{2~O_{cb}~O_{ba}O_{ca}}.          \eqno(5.13)
$$
Now in the same way let us calculate the next ratio, $ \frac{\Pi_{bc}}{O_{bc}}.$ Formula for the projection  evaluated
as follows
$$
\PPI_{bc}=r_a\cos a_1-r_a\cos a_2=
$$
$$
\frac{O^2_{ba}-r^2_b+r^2_a}{2~O_{ba}}-\frac{O^2_{ca}-r^2_c+r^2_a}{2O_{ac}}=
$$
$$
\frac{~(~O^2_{ba}-r^2_b+r^2_a~)~O_{ac} -(~O^2_{ca}-r^2_c+r^2_a~)O_{ba}}{2O_{ac}O_{ba}}=
$$
$$
\frac{~(~O^2_{ba}O_{ac}-r^2_bO_{ac}+r^2_a~O_{ac}~) -(~O^2_{ca}O_{ba}-r^2_cO_{ba}+r^2_a~O_{ba}~)}{2O_{ac}O_{ba}}=
$$
$$
\frac{~(~O_{ba}O_{ac}~(O_{ba}-O_{ca}) -r^2_bO_{ac}+r^2_a~O_{ac}~) +r^2_cO_{ba}-r^2_a~O_{ba}~)}{2O_{ac}O_{ba}}.
$$
Take into account $-O_{cb}=O_{ba}-O_{ac}$, hence,
$$
\PPI_{bc}  \frac{~-O_{ba}O_{ac}~(O_{cb} -r^2_bO_{ac}+r^2_a~(O_{cb})+r^2_cO_{ba}~)}{2O_{ac}O_{ba}}=
$$
$$
\frac{~(~O_{ba}O_{ac}~(O_{cb}) -r^2_bO_{ac}+r^2_a~(O_{cb})+r^2_cO_{ba}~)}{2O_{ac}O_{ba}}.
$$
Now, calculate the ration
$$
\frac{\PPI_{bc}}{O_{bc}}= \frac{~(~O_{ba}O_{ac}~(O_{cb})
-r^2_bO_{ac}+r^2_a~(O_{cb})+r^2_cO_{ba}~)}{2O_{ac}O_{ba}O_{bc}}.\eqno(5.14)
$$
This expression coincides with (5.13), consequently,
$$
\frac{\PPI_{bc}}{O_{bc}}= \frac{\PPI_{ca}}{O_{ca}}. \eqno(5.15)
$$
By taking into account (5.7), we arrive to the desired relations
$$
\frac{\PPI(BC) }{ O_{bc}}=\frac{\PPI(AC)}{O_{ca}}=\frac{\PPI(AB)}{O_{ab}}. \eqno(5.16)
$$

{\bf End of proof.}

Now come back to designations introduced in Section 3. In these designations equations (5.16) are written as follows
$$
\frac{w_{01}-v_{01}}{~w_{03}-w_{02}~}= \frac{w_{02}-v_{02}}{~v_{03}+v_{01}~}=
\frac{w_{03}-v_{03}}{~w_{01}+v_{02}~}.\eqno(5.17)
$$

{\bf Theorem 5.2}

{\it The sides and the angles  of triangle $\widetilde{\triangle} ABC$ satisfy the equations (5.1).}

{\bf Proof}

By using the designations introduced in Sec.2 the system of equations (6.1) can be written as follows
$$
\frac{\sinh a}{\sin\alpha}= \frac{r_a(w_{01}-v_{01})}{~yz}~:~\frac{x~(~w_{03}-w_{02}~)}{r_br_c}=
\frac{w_{01}-v_{01}}{~w_{03}-w_{02}~}\frac{xyz}{r_ar_br_c},
$$
$$
\frac{\sinh b}{\sin\beta}= \frac{w_{02}-v_{02}}{r_b~xz}~:~\frac{y~(~v_{03}+v_{01}~)}{r_ar_c}=
\frac{w_{02}-v_{02}}{~v_{03}+v_{01}~}\frac{xyz}{r_ar_br_c},
$$
$$
\frac{\sinh c}{\sin\delta}= \frac{w_{03}-v_{03}}{r_c~yx}~:~\frac{z~(~w_{01}+v_{02}~)}{r_br_a}=
\frac{w_{03}-v_{03}}{~w_{01}+v_{02}~}\frac{xyz}{r_ar_br_c}. \eqno(5.18)
$$
It is seen, these equations contain a common factor which is symmetric with respect to $a,b,c$ and $x,y,z$. Multiply
all equations (5.17) by this factor. We arrive to equations (5.1).

{\bf End of proof}

{\bf Theorem 5.3}.

{\it The sides and the angles  of triangle $\widetilde{\triangle} ABC$ satisfy the following equation.}
$$
\cos\delta=\sin\alpha\sin\beta\cosh c-\cos\alpha\cos\beta. \eqno(5.19)
$$

{\bf Proof.}

Evaluate the first term of the right-hand side of (5.19).
$$
\sin\alpha\sin\beta\cosh c=\frac{1}{r_ar_b}O_{ca}~O_{cb}~(1-\cos c_1\cos c_2)=
$$
$$
\frac{1}{r_ar_b}O_{ca}~O_{cb}(1-\frac{O^2_{cb}-r^2_b+r^2_c}{2r_cO_{cb}}~~\frac{O^2_{ca}-r^2_a+r^2_c}{2r_c~O_{ac}} ) =
$$
$$
\frac{1}{r_ar_b}O_{ca}~O_{cb}-\frac{1}{r_ar_b}O_{ca}~O_{cb}
\frac{O^2_{cb}-r^2_b+r^2_c}{2r_cO_{cb}}~~\frac{O^2_{ca}-r^2_a+r^2_c}{2r_c~O_{ac}} ) =
$$
$$
\frac{1}{r_ar_b}O_{ca}~O_{cb}-\frac{1}{4r_ar_br_c^2} (~{O^2_{cb}-r^2_b+r^2_c}~)~(~{O^2_{ca}-r^2_a+r^2_c}~) =
$$
$$
\frac{1}{4r_ar_br_c^2}~(~~4O_{ca}~O_{cb}~r_c^2-  (~{O^2_{cb}-r^2_b+r^2_c}~)~(~{O^2_{ca}-r^2_a+r^2_c}~)~)=
$$
$$
\frac{1}{4r_ar_br_c^2}~(~~4O_{ca}~O_{cb}~r_c^2-~O^2_{cb}O^2_{ca}-(~O^2_{cb}+O^2_{ca})r^2_c+O^2_{cb}r^2_a~+O^2_{ca}r^2_b).
\eqno(5.20)
$$

Now calculate the product  $ \cos\alpha\cos\beta$ by using the following formulae
$$
\cos\alpha=\frac{1}{2r_br_c}(r_c^2+r_b^2-O^2_{cb}),~~~ \cos\beta=-\frac{1}{2r_ar_c}(r_c^2+r_a^2-O^2_{ca}).
$$
We get
$$
\cos\alpha\cos\beta=\frac{1}{2r_br_c}(r_c^2+r_b^2-O^2_{cb})~\frac{1}{2r_ar_c}(r_c^2+r_a^2-O^2_{ca})=
$$
$$
\frac{1}{4r_br_ar^2_c}( ~~O^2_{cb}O^2_{ca}~-(~O^2_{cb}+O^2_{ca})r^2_c-O^2_{cb}r^2_a-O^2_{ca}r^2_b). \eqno(5.21)
$$
By using equations (5.20) and (5.21) calculate the difference
$$
\sin\alpha\sin\beta\cosh c-\cos\alpha\cos\beta=
$$
$$
\frac{1}{4r_ar_br_c^2}~(~~4O_{ca}~O_{cb}~r_c^2-~O^2_{cb}O^2_{ca}-(~O^2_{cb}+O^2_{ca})r^2_c+O^2_{cb}r^2_a~+O^2_{ca}r^2_b)-
$$
$$
-\frac{1}{4r_br_ar^2_c}( ~+O^2_{cb}O^2_{ca}~-(~O^2_{cb}+O^2_{ca})r^2_c-O^2_{cb}r^2_a-O^2_{ca}r^2_b)
$$
$$
=\frac{1}{2r_br_a}(r_a^2+r_b^2-O^2_{ab})
$$
$$
=\cos\delta.  \eqno(5.22)
$$
Thus,  we got the equation (5.19).

{\bf End of proof}.

The other two equations, (5.3)and (5.4), are proved analogously.

{\bf Concluding remarks.}

We have seen that the hyperbolic trigonometry, like circular angle, gives arise in a natural way on the Euclidean
plane. The hyperbolic description of the elements of the Euclidean plane has to be considered as a complementary tool
of the Euclidean Geometry. This description provides with new insights into hidden nature of the Euclidean Geometry.

The proofs of theorems "hyperbolic law of cosines I", "hyperbolic law of sines" and "hyperbolic law od cosines II" were
based purely on elements of the Euclidean geometry. These laws express interrelations between distances of the circles,
radii and angles between radiuses and $X$-axis.

The method developed in this paper opens new pathway from Euclidean to hyperbolic geometry and can be used as an
introduction into complex field of hyperbolic geometry.

\section{Appendix }

The task of this section is to reduce the expression
$$
{v_2w_1}~(~w^2_{2}+v^2_{1}+(w_{02}-v_{02})^2+(w_{01}-v_{01})^2~)-2v_{02}w_{01}(w_{02}-v_{02})(w_{01}-v_{01}),\eqno(A.1)
$$
onto the expression
$$
2(r_a^2-v_{01}w_{01})(~r_b^2 -v_{02}w_{02}~). \eqno(A.2)
$$

By taking into account the equation $ v_2=w_1 $ and opening the brackets transform (A.1) into the following form
$$
=\overline{w^2_1w_2^2}+\widetilde{v_2^2v_1^2}+
$$
$$
\overline{w_1^2~w^2_{02}}+w_1^2~v^2_{02}-\underline{2~w_{02}v_{02}~w_1^2}+
$$
$$
w^2_{01}~v_2^2+\widetilde{v^2_{01}v_2^2}-\underbrace{2v_{01}w_{01}~v_2^2}+
$$
$$
2v_{02}w_{01}w_{02}v_{01}-
$$
$$
\underline{2v_{02}w_{02}~w^2_{01}}+2v^2_{02}w^2_{01}-\underbrace{2v^2_{02}w_{01}v_{01}}.\eqno(A.3)
$$
Joint together terms marked by same under- and over- lines and take into account that
$$
\overline{w^2_1w_2^2}+\overline{w_1^2~w^2_{02}}=w_1^2~r_b^2
$$
$$
 \widetilde{v_2^2v_1^2}+\widetilde{v^2_{01}v_2^2}=v_2^2~r_a^2
$$
$$
-\underline{2~w_{02}v_{02}~w_1^2}-\underline{2v_{02}w_{02}~w^2_{01}}=-2v_{02}w_{02}~r_a^2
$$
$$
-\underbrace{2v_{01}w_{01}~v_2^2}-\underbrace{2v^2_{02}w_{01}v_{01}}=-2v_{01}w_{01}~r_b^2.
$$\\

The last term in (A.3) represent as follows
$$
(2v^2_{02}w^2_{01}={v^2_{02}w^2_{01}}+{v^2_{02}w^2_{01}}~).
$$
In this way we transform expression (A.3) into the following form
$$
2v_{02}w_{01}w_{02}v_{01}
$$
$$
-2v_{02}w_{02}~r_a^2 -2v_{01}w_{01}~r_b^2+
$$
$$
w_1^2~r_b^2+v_2^2~r_a^2+
$$
$$
\underbrace{w_1^2~v^2_{02}}+\overbrace{w^2_{01}~v_2^2}+
$$
$$
\overbrace{v^2_{02}w^2_{01}}+\underbrace{v^2_{02}w^2_{01}}~.\eqno(A.4)
$$\\
Join term with same under- and over- lines taking into account that
$$
\overbrace{w^2_{01}~v_2^2}+ \overbrace{v^2_{02}w^2_{01}}= w^2_{01}~r_b^2,
$$
$$
\underbrace{w_1^2~v^2_{02}}+\underbrace{v^2_{02}w^2_{01}}= v^2_{02}~r_a^2.
$$
Then, we fulfil to the following set of simple transformations
$$
2v_{02}w_{01}w_{02}v_{01}
$$
$$
-2v_{02}w_{02}~r_a^2 -2v_{01}w_{01}~r_b^2
$$
$$
+ w_1^2~r_b^2+v_2^2~r_a^2
$$
$$
 + v^2_{02}~r_a^2+ w^2_{01}~r_b^2=
$$
$$
2v_{02}w_{01}w_{02}v_{01}+
$$
$$
-2v_{02}w_{02}~r_a^2 -2v_{01}w_{01}~r_b^2
$$
$$
2r_a^2~r_b^2=
$$
$$
2(r_a^2-v_{01}w_{01})(~r_b^2 -v_{02}w_{02}~)
$$

{\bf End of proof}.

\end{document}